\documentclass[11pt]{article}

\usepackage[margin=1in]{geometry}
\usepackage{amsmath,amssymb,amsthm,mathtools}
\usepackage{enumitem}
\usepackage{booktabs}
\usepackage{tikz}
\usetikzlibrary{arrows.meta,calc,positioning}
\usepackage{pgfplots}
\pgfplotsset{compat=1.17}
\usepackage{hyperref}

\hypersetup{
  colorlinks=true,
  linkcolor=blue,
  citecolor=blue,
  urlcolor=blue
}

\newtheorem{theorem}{Theorem}[section]
\newtheorem{lemma}[theorem]{Lemma}
\newtheorem{proposition}[theorem]{Proposition}
\newtheorem{corollary}[theorem]{Corollary}
\theoremstyle{definition}
\newtheorem{definition}[theorem]{Definition}

\theoremstyle{remark}
\newtheorem{remark}[theorem]{Remark}

\newcommand{\Id}{\mathrm{Id}}
\newcommand{\eps}{\varepsilon}

\newcommand{\J}{\mathcal{J}}
\newcommand{\R}{\mathbb{R}}

\title{Stability in the Non-Neutral and Neutral Regimes of the Holonomy-Based Inverse Problem in the Class \(\mathcal{O}_C\)}
\author{M. Barkatou\thanks{ISTM Laboratory, Chouaib Doukkali University, Morocco. \texttt{barkatou.m@ucd.ac.ma}} \and M. El Morsalani\thanks{QWAVE-CONSULT.EU. \texttt{Mohamed.elmorsalani@qwave-consult.eu}}}
\date{\today}

\begin{document}

\maketitle

\begin{abstract}
We study stability of anchor reconstruction in the holonomy formulation of the inverse problem for return maps in the class $\mathcal O_C$. For a known convex core $C$ and an observed return map, reconstruction of the unknown thickness reduces to recovery of an anchor value $a_*$ that is a fixed point of nonlinear loop transports $T_\gamma$. Writing $H=T_\gamma-\Id$, local identifiability is organized by the first nonzero derivative of $H$ at $a_*$.

In the non-neutral case $H'(a_*)\ne0$, we prove local Lipschitz stability and an observable multiplier estimate. In the neutral regime, we derive the exact second and third variations of the holonomy map and, more generally, the full finite variational jet through the Bell-polynomial form of the Fa\`a di Bruno equations. If $H(a_*)=0$, $H^{(j)}(a_*)=0$ for $1\le j<k$, and $H^{(k)}(a_*)\ne0$, every sufficiently localized perturbed anchor satisfies a quantitative $1/k$-H\"older estimate; odd and even orders have different existence mechanisms. We also prove stability of every fixed finite holonomy jet under perturbations of the observed return data and deterministic multi-loop aggregation at the minimal isolation order. For finite-dimensional $C^\infty$ families of admissible return data, Malgrange preparation reduces the local perturbed anchor equation at an order-$k$ holonomy singularity to a smooth nonvanishing factor times a monic polynomial of degree $k$ in the anchor displacement, exposing the discriminant geometry underlying root splitting. We further prove real-analytic dependence of the loop transport on its initial anchor away from the singular fibre set. Consequently, infinite-order flatness for one loop implies that the corresponding holonomy map is locally the identity.

We also prove that exact uniform separation persists under sufficiently small return-data perturbations. Under this quantitative admissibility condition, the local anchor estimate propagates to the full reconstructed thickness, yielding global $C^0$-stability of order $O(\varepsilon_F^{1/k})$ at finite order. The theory is validated numerically on exact scalar benchmarks and on two genuine circular-core return geometries: a non-neutral example exhibiting first-order convergence and a numerically tuned neutral example with $H'(a_*)$ at machine zero and $H''(a_*)\ne0$, for which both anchor and global thickness errors exhibit the predicted square-root law.
\end{abstract}



\tableofcontents

\section{Introduction}
\label{sec:introduction}

The return-map construction in the class $\mathcal O_C$ associates to a known convex core $C$ and a surrounding admissible normal-graph boundary a discrete map $F:\partial C\to\partial C$. The forward geometry and regularity of this map were developed in \cite{Barkatou2002,BarkatouElMorsalaniReturn2026}. The inverse problem studied in \cite{BarkatouElMorsalaniInv2026} asks how much of the unknown thickness function $d$, and hence of the outer boundary, can be reconstructed when $C$ and $F$ are known.

A central result of that inverse theory is that exact non-uniqueness reduces to a scalar anchor value. The inverse system transports a trial value $a=d(c_0)$ along paths of $M:=\partial C$, and closed paths generate nonlinear holonomy maps $T_\gamma$. The true anchor $a_*$ is a common fixed point. When
\[
H(a):=T_\gamma(a)-a
\]
satisfies $H'(a_*)\ne0$, the anchor is locally isolated by the ordinary transverse-root mechanism. The singular regime begins when $H'(a_*)=0$: linearized holonomy becomes neutral, but higher vertical variations of the nonlinear transport may still retain identifying information.

The word ``holonomy'' also appears in several nearby geometric inverse problems, but the data and unknowns are different. Ceki\'c and Lefeuvre recover unitary connections, modulo gauge, from traces of connection holonomies along primitive closed geodesics and prove local injectivity near generic connections \cite{CekicLefeuvre2025}; their companion stability paper establishes H\"older-type estimates for that problem \cite{CekicLefeuvre2024}. More broadly, transport equations, scattering data, non-Abelian X-ray transforms, boundary rigidity and range questions form a central part of contemporary geometric inverse problems \cite{PaternainSaloUhlmann2023,Bohr2021,CekicGuillarmouLefeuvre2024,Grebnev2026}. In the present paper, by contrast, the holonomy is a scalar nonlinear transport in the fibre variable representing the unknown thickness anchor. This distinction is important: the H\"older exponents below arise from finite-order singularity of a scalar anchor fixed point, while their coefficients are generated by the return-map transport geometry.

\paragraph{Contribution and positioning.}
The $1/k$ exponent itself should not be confused with the novelty of the paper. Multiple roots are classically ill-conditioned and their numerical treatment requires multiplicity-aware analysis; see, for example, \cite{Zeng2005}. The purpose here is to close, within the return-map inverse problem, the problem-specific chain
\[
\begin{aligned}
\text{return-data perturbation}
&\Longrightarrow \text{finite holonomy-jet perturbation}
\Longrightarrow \text{isolation order},\\
&\Longrightarrow \text{anchor stability}
\Longrightarrow \text{boundary reconstruction stability}.
\end{aligned}
\]
The genuinely new ingredients are the geometric realization of the multiplicity through the exact holonomy jet, its finite-jet stability under noisy return data, the Malgrange-prepared local polynomial unfolding for finite-dimensional smooth families, deterministic multi-loop aggregation, and the propagation of the anchor estimate to the full reconstructed thickness. The scalar multiple-root exponent is classical; its geometric coefficient, observability, and inverse-problem consequences are the contributions developed below.

Define
\[
k_\gamma(a_*):=\min\{k\ge1:H^{(k)}(a_*)\ne0\},
\]
with $k_\gamma(a_*)=+\infty$ if all derivatives vanish. The main results are as follows.
\begin{enumerate}[label=(\roman*)]
\item In the non-neutral case we prove a local one-loop stability theorem, an observable multiplier bound, and quantitative transfer of return-map errors to holonomy errors.
\item We derive the exact second holonomy variation and the explicit third variation, and prove the full finite variational hierarchy
\[
Z_n'=\sum_{r=1}^n G_{y^r}\,B_{n,r}(Z_1,\ldots,Z_{n-r+1}),
\]
where $B_{n,r}$ are partial Bell polynomials.
\item If the first nonzero holonomy derivative has order $k$, every sufficiently localized perturbed zero satisfies a quantitative $1/k$-H\"older estimate. Odd and even orders have different existence mechanisms.
\item For every finite-dimensional $C^\infty$ family of admissible return data passing through an order-$k$ holonomy singularity, the Malgrange preparation theorem gives a local factorization into a smooth unit and a monic degree-$k$ polynomial in the anchor displacement. Thus the perturbed anchor set is locally governed by a finite polynomial unfolding, while our stability estimates quantify the conditioning of its roots.
\item For every fixed finite $p$, $C^0$ perturbation of the observable return components controls the holonomy jet $(T_\gamma,T_\gamma',\ldots,T_\gamma^{(p)})$ on uniformly separated anchor intervals. This makes the higher-order stability coefficients observable from noisy return data.
\item We prove a deterministic multi-loop estimate at the minimal isolation order, including approximate common fixed points.
\item We prove real-analytic dependence of $T_\gamma$ on the initial anchor on uniformly separated intervals. Thus $k_\gamma(a_*)=+\infty$ implies that $T_\gamma$ is locally the identity near $a_*$.
\item Under uniform separation from the singular fibre set, the local anchor estimate propagates to a global estimate
\[
\|\hat d-d\|_{C^0(M)}=O(\varepsilon_F^{1/k})
\]
for every fixed finite order $k\ge2$.
\item We validate both the non-neutral and neutral regimes on genuine return geometries of the unit circular core. The observed convergence rates are approximately $1.0002$ in the non-neutral case and $0.4960$ in the tuned neutral case; the corresponding global thickness rates are $1.0002$ and $0.4953$.
\end{enumerate}

The geometric origin of the higher-order coefficients beyond the explicit vertical variational formulas remains subtler. In particular, invariant small-loop formulas on a general non-umbilic core must retain the Levi-Civita connection terms of the principal frame. We therefore separate the proved holonomy-jet theory from those remaining geometric genericity questions.

The paper is organized as follows. Section~\ref{sec:framework} recalls the nonlinear transport formulation. Section~\ref{sec:nonneutral} treats the non-neutral regime. Section~\ref{sec:variational} develops the variational holonomy jet, its Malgrange-prepared finite-dimensional unfoldings, and the analytic flat regime. Section~\ref{sec:finiteorder} proves finite-order isolation and parity-sensitive existence. Section~\ref{sec:observable-k} transfers return-data perturbations to finite holonomy jets and develops observable coefficients. Section~\ref{sec:multiloop} gives deterministic multi-loop aggregation. Section~\ref{sec:globalneutral} propagates anchor errors globally. Section~\ref{sec:geometry} records the remaining invariant geometric questions. Section~\ref{sec:numerics} contains executed algebraic and geometric numerical experiments, including a TikZ visualization of the tuned neutral geometry.

\section{Framework and nonlinear transport}
\label{sec:framework}

We use the notation and standing hypotheses of the inverse problem developed in \cite{BarkatouElMorsalaniInv2026}, within the class \(\mathcal O_C\) introduced in \cite{Barkatou2002}. Let \(C\subset\R^N\) be compact convex with \(C^3\) boundary \(M:=\partial C\). Let \(d:M\to(0,\infty)\) be of class \(C^2\), let \(F:M\to M\) be the associated return map, and let \(S_c:=-D\nu(c)|_{T_cM}\) be the shape operator of \(M\). We denote by \(u_F,q_F\) the tangential and normal components of \(F(c)-c\), and by \(T_\gamma\) the holonomy transport along a closed \(C^1\) loop \(\gamma\) based at \(c_0\).

For a fixed closed loop \(\gamma:[0,1]\to M\) based at \(c_0\), we write the scalar transport equation as
\begin{equation}
y_a'(t)=G_\gamma(t,y_a(t)),\qquad y_a(0)=a,
\end{equation}
where
\begin{equation}
G_\gamma(t,y)=\frac{A(t)-yB(t)}{y-q(t)},
\end{equation}
with
\begin{equation}
A(t)=\langle u_F(\gamma(t)),\dot\gamma(t)\rangle,\qquad
B(t)=\langle S_{\gamma(t)}u_F(\gamma(t)),\dot\gamma(t)\rangle,\qquad
q(t)=q_F(\gamma(t)).
\end{equation}
Set
\begin{equation}
P(t)=A(t)-q(t)B(t)
=\langle (I-q(t)S_{\gamma(t)})u_F(\gamma(t)),\dot\gamma(t)\rangle.
\end{equation}
On the domain \(y-q(t)>0\),
\begin{equation}
\partial_yG_\gamma(t,y)=-\frac{P(t)}{(y-q(t))^2},
\end{equation}
\begin{equation}
\partial_y^2G_\gamma(t,y)=2\frac{P(t)}{(y-q(t))^3},
\end{equation}
and more generally, for \(r\ge 1\),
\begin{equation}
\partial_y^rG_\gamma(t,y)=(-1)^r r!\frac{P(t)}{(y-q(t))^{r+1}}.
\end{equation}

The holonomy displacement is
\begin{equation}
H(a):=T_\gamma(a)-a,\qquad T_\gamma(a)=y_a(1),
\end{equation}
with \(H(a_*)=0\) at the true anchor \(a_*:=d(c_0)\). We set
\begin{equation}
\delta_\gamma:=|H'(a_*)|.
\end{equation}
The isolation order is
\begin{equation}
k_\gamma(a_*):=\min\{k\ge 1:H^{(k)}(a_*)\neq 0\},
\end{equation}
with \(k_\gamma(a_*)=+\infty\) if all derivatives vanish. When \(k<\infty\), we write
\begin{equation}
m_k:=|H^{(k)}(a_*)|>0.
\end{equation}

We use two error norms on the holonomy displacement:
\begin{equation}
\eps:=\|\hat H-H\|_{C^0(\bar I)},\qquad
\eps_T:=\|\hat H-H\|_{C^1(\bar I)}.
\end{equation}
Evidently \(\eps\le\eps_T\). The non-neutral regime requires \(C^1\) control and is stated with \(\eps_T\). The neutral regime requires only \(C^0\) control plus a Taylor remainder, and is stated with \(\eps\); \(C^1\) control is invoked only when existence of the perturbed fixed point is established by a contraction.

For \(j\ge 1\), we set
\begin{equation}
M_j:=\|H^{(j)}\|_{C^0(\bar I)}.
\end{equation}

\subsection{Uniform separation and persistence under perturbation}
\label{subsec:separation}
All stability statements below require the horizontal lifts to remain at positive distance from the singular fibre set $y=q_F(\gamma(t))$. We first record that this quantitative admissibility condition is stable under sufficiently small $C^0$ perturbations of the observed return components.

\begin{definition}[Uniform separation]
\label{def:uniform-separation}
Let $I_0\subset\mathbb R$ be a compact anchor interval and $\gamma:[0,1]\to M$ a $C^1$ curve. We say that the exact transport is uniformly separated on $I_0$ with margin $\eta>0$ if, for every $a\in I_0$, the lift $y_a$ exists on $[0,1]$ and
\[
y_a(t)-q_F(\gamma(t))\ge\eta
\qquad\text{for all }(a,t)\in I_0\times[0,1].
\]
\end{definition}

\begin{lemma}[Persistence of uniform separation]
\label{lem:separation-persistence}
Assume Definition~\ref{def:uniform-separation} holds with margin $\eta>0$. Set
\[
\varepsilon_F:=\|\hat u_F-u_F\|_{C^0(\gamma)}+\|\hat q_F-q_F\|_{C^0(\gamma)}.
\]
Then there exists $\varepsilon_F^*>0$, depending only on $\eta$, $I_0$, $\gamma$ and uniform bounds for the exact transport on the compact family of exact solution graphs, such that, whenever $0\le\varepsilon_F\le\varepsilon_F^*$, every perturbed lift $\hat y_a$ with $a\in I_0$ exists on $[0,1]$ and
\[
\hat y_a(t)-\hat q_F(\gamma(t))\ge\frac{\eta}{2}
\qquad\text{for all }(a,t)\in I_0\times[0,1].
\]
Moreover,
\[
\sup_{a\in I_0}\sup_{0\le t\le1}|\hat y_a(t)-y_a(t)|\le C_{\gamma,I_0}\varepsilon_F
\]
for a constant independent of the perturbation once $\varepsilon_F\le\varepsilon_F^*$.
\end{lemma}

\begin{proof}
Because $I_0$ is compact and the exact lifts depend continuously on $(a,t)$, their union is contained in a compact fibre interval. Fix the tube
\[
\mathcal T:=\{(a,t,y):a\in I_0,\ 0\le t\le1,\ |y-y_a(t)|\le\eta/4\}.
\]
If $\varepsilon_F\le\eta/4$, then throughout $\mathcal T$ one has
\[
y-\hat q_F(\gamma(t))
\ge y_a(t)-q_F(\gamma(t))-\frac{\eta}{4}-\varepsilon_F
\ge\frac{\eta}{2}.
\]
Hence the exact and perturbed rational vector fields are uniformly bounded and Lipschitz in $y$ on this tube, and there is a constant $C_0$ such that
\[
|\hat G(t,y)-G(t,y)|\le C_0|\dot\gamma(t)|\varepsilon_F
\qquad\text{on }\mathcal T.
\]
For fixed $a\in I_0$, let $\tau$ be the supremum of times for which the perturbed solution exists and remains in the tube, equivalently $|\hat y_a-y_a|\le\eta/4$, on $[0,\tau]$. On this interval the standard Gronwall comparison gives
\[
|\hat y_a(t)-y_a(t)|
\le C_0\ell(\gamma)e^{L_0\ell(\gamma)}\varepsilon_F,
\qquad 0\le t\le\tau,
\]
with a uniform Lipschitz constant $L_0$ on $\mathcal T$. Choose $\varepsilon_F^*$ so that the right-hand side is strictly smaller than $\eta/4$ and also $\varepsilon_F^*\le\eta/4$. A first exit at $\tau<1$ is then impossible. Thus $\tau=1$, the displayed difference estimate holds on the full interval, and
\[
\hat y_a(t)-\hat q_F(\gamma(t))
\ge\eta-\frac{\eta}{4}-\frac{\eta}{4}
=\frac{\eta}{2}.
\]
All constants are uniform in $a\in I_0$, which completes the proof.
\end{proof}

\begin{remark}[Conditioning near the singular fibre]
\label{rem:separation-conditioning}
The vertical Lipschitz constants below scale like $\eta^{-2}$. Consequently the Gronwall amplification factors deteriorate rapidly as the exact separation margin tends to zero. This blow-up records genuine ill-conditioning near the non-admissible set $y=q_F$, rather than a proof artefact.
\end{remark}

\begin{remark}[Certifying separation near the true anchor]
For a concrete geometry it is enough first to certify a positive separation margin along the true horizontal section $y_*(t)=d(\gamma(t))$. Continuity of the transport with respect to the initial anchor then gives a compact interval $I_0$ around $a_*$ on which Definition~\ref{def:uniform-separation} holds, possibly with a smaller margin. In Section~\ref{sec:numerics}, the interval enclosures of the circular-core geometry provide precisely such quantitative geometric margins; they are used as evidence for a well-separated true section, not as a substitute for Definition~\ref{def:uniform-separation} on an arbitrary anchor interval.
\end{remark}

\section{Non-neutral stability}
\label{sec:nonneutral}

Throughout this section we assume
\[
\delta_\gamma:=|H'(a_*)|>0.
\]
This is the non-neutral regime already underlying the first-order holonomy mechanism in \cite{BarkatouElMorsalaniInv2026}.  The purpose here is quantitative stability under perturbed return data.

\begin{theorem}[Local one-loop stability in the non-neutral regime]
\label{thm:nonneutral}
Let \(H=T_\gamma-\Id\) and \(\hat H=\hat T_\gamma-\Id\) on an interval \(I\). Assume
\[
H(a_*)=0,\qquad \delta_\gamma:=|H'(a_*)|>0,
\]
\[
\|\hat H-H\|_{C^1(\bar I)}\le \eps_T,
\qquad
\|H''\|_{C^0(\bar I)}\le M_2,
\]
and
\[
6M_2\eps_T<\delta_\gamma^2,
\qquad
4\eps_T<\delta_\gamma.
\]
Set
\[
r:=\frac{3\eps_T}{\delta_\gamma}
\]
and assume
\[
\overline{B_r(a_*)}\subset I.
\]
Then \(\hat H\) has a unique zero \(\hat a\) in \(\overline{B_r(a_*)}\); in fact \(\hat a\in B_r(a_*)\), and
\begin{equation}
|\hat a-a_*|
\le
\frac{\eps_T}{\delta_\gamma-\eps_T-3M_2\eps_T/\delta_\gamma}.
\label{eq:nonneutral-refined}
\end{equation}
In particular,
\begin{equation}
|\hat a-a_*|\le \frac{4\eps_T}{\delta_\gamma}.
\label{eq:nonneutral-simple}
\end{equation}
If \(\eps_T/\delta_\gamma\ll1\) and \(M_2\eps_T/\delta_\gamma^2\ll1\), then only the upper asymptotic estimate
\[
|\hat a-a_*|
\le
\frac{\eps_T}{\delta_\gamma}
\left(1+O\left(\frac{\eps_T}{\delta_\gamma}
+\frac{M_2\eps_T}{\delta_\gamma^2}\right)\right)
\]
is asserted. No global uniqueness in \(I\) follows from these hypotheses alone.
\end{theorem}

\begin{proof}
Set \(r=3\eps_T/\delta_\gamma\). For \(a\in\overline{B_r(a_*)}\),
\[
|H'(a)|\ge\delta_\gamma-M_2r>\frac{\delta_\gamma}{2},
\]
so \(H'\) has constant sign on the closed interval \(\overline{B_r(a_*)}\). Replacing \(H\) and \(\hat H\) by their negatives if necessary, suppose \(H'>0\).
Moreover,
\[
H(a_*+r)
=\int_0^rH'(a_*+s)\,ds
\ge \delta_\gamma r-\frac12M_2r^2
>\eps_T,
\]
and similarly \(H(a_*-r)<-\eps_T\). Hence \([-\eps_T,\eps_T]\subset H(\overline{B_r(a_*)})\), and the inverse branch \(H^{-1}\) is defined there with
\[
|(H^{-1})'(e)|\le\frac{2}{\delta_\gamma}.
\]
Define on the complete metric space \(\overline{B_r(a_*)}\)
\[
\Phi(a):=H^{-1}\bigl(-(\hat H-H)(a)\bigr).
\]
Then
\[
|\Phi(a)-a_*|
\le\frac{2\eps_T}{\delta_\gamma}
<r,
\]
so \(\Phi(\overline{B_r(a_*)})\subset B_r(a_*)\). Also
\[
|\Phi'(a)|
\le\frac{2\eps_T}{\delta_\gamma}<\frac12.
\]
Banach's fixed-point theorem therefore gives a unique fixed point \(\hat a\in\overline{B_r(a_*)}\), necessarily in the open ball. The identity \(\Phi(\hat a)=\hat a\) is equivalent to \(\hat H(\hat a)=0\), and every zero of \(\hat H\) in the closed ball is a fixed point of \(\Phi\), proving uniqueness.

Finally,
\[
0=\hat H(\hat a)
=\hat H(a_*)+\int_{a_*}^{\hat a}\hat H'(s)\,ds.
\]
On the segment joining \(a_*\) to \(\hat a\),
\[
|\hat H'(s)|\ge\delta_\gamma-M_2r-\eps_T,
\]
while \(|\hat H(a_*)|\le\eps_T\). This yields \eqref{eq:nonneutral-refined}; the remaining estimates follow immediately.
\end{proof}

\subsection{Observable first-order coefficient}

The exact quantity \(\delta_\gamma\) depends on the unknown anchor. Define instead
\[
\hat\delta_\gamma:=|\hat T_\gamma'(\hat a)-1|=|\hat H'(\hat a)|.
\]

\begin{lemma}[Proximity of the observed multiplier]
\label{lem:proximity}
Under the hypotheses of Theorem~\ref{thm:nonneutral},
\[
|\hat\delta_\gamma-\delta_\gamma|
\le
\eps_T+M_2|\hat a-a_*|
\le
\eps_T\left(1+\frac{3M_2}{\delta_\gamma}\right).
\]
\end{lemma}

\begin{proof}
Use the reverse triangle inequality and split
\[
\hat H'(\hat a)-H'(a_*)
=[\hat H'(\hat a)-H'(\hat a)]+[H'(\hat a)-H'(a_*)].
\]
The first term is bounded by \(\eps_T\), the second by \(M_2|\hat a-a_*|\), and \(|\hat a-a_*|<r\).
\end{proof}

\begin{corollary}[Observable non-neutral bound]
\label{cor:observable-nonneutral}
Assume, in addition,
\[
\eps_T\left(1+\frac{3M_2}{\delta_\gamma}\right)
\le\frac{\delta_\gamma}{2}.
\]
Then
\[
\frac{\delta_\gamma}{2}
\le\hat\delta_\gamma
\le\frac{3\delta_\gamma}{2},
\]
and the anchor error admits the observable-form estimate
\begin{equation}
|\hat a-a_*|
\le
\frac{9}{2}\frac{\eps_T}{\hat\delta_\gamma}.
\label{eq:observable-delta}
\end{equation}
\end{corollary}

\begin{proof}
The interval for \(\hat\delta_\gamma\) follows from Lemma~\ref{lem:proximity}. Since \(\hat\delta_\gamma\le 3\delta_\gamma/2\), one has \(\delta_\gamma\ge2\hat\delta_\gamma/3\). Combining this with \(|\hat a-a_*|<3\eps_T/\delta_\gamma\) gives \eqref{eq:observable-delta}.
\end{proof}

\subsection{From return-map perturbations to holonomy perturbations}

\begin{lemma}[Continuity of the transport map]
\label{lem:reduction}
Let \(\gamma:[0,1]\to M\) be a \(C^1\) curve, and set
\[
\eps_F:=\|\hat u_F-u_F\|_{C^0(\gamma)}
+\|\hat q_F-q_F\|_{C^0(\gamma)}.
\]
Fix \(\eps_0>0\) and assume \(0\le\eps_F\le\eps_0\). Assume that the exact transport issued from the compact interval \(I_0\) is uniformly separated with margin \(\eta>0\) in the sense of Definition~\ref{def:uniform-separation}. Shrinking \(\eps_0\) if necessary, Lemma~\ref{lem:separation-persistence} gives perturbed separation with margin \(\eta/2\) and a common bounded fibre region. Then there exist constants \(C_0,C_1,L_0,L_1\), depending only on that common fibre bound, \(\eta\), \(\eps_0\), the geometry of \(M\), and uniform bounds for the exact fields, such that along \(\gamma\)
\[
|\partial_yG(t,y)|+|\partial_y\hat G(t,y)|
\le L_0|\dot\gamma(t)|,
\]
\[
|\partial_y^2\hat G(t,y)|\le L_1|\dot\gamma(t)|,
\]
\[
|\hat G(t,y)-G(t,y)|\le C_0|\dot\gamma(t)|\eps_F,
\]
and
\[
|\partial_y\hat G(t,y)-\partial_yG(t,y)|
\le C_1|\dot\gamma(t)|\eps_F.
\]
Consequently, with \(\ell(\gamma)=\int_0^1|\dot\gamma|dt\),
\begin{align}
\|\hat T_\gamma-T_\gamma\|_{C^0(I_0)}
&\le C_0\ell(\gamma)e^{L_0\ell(\gamma)}\eps_F,\label{eq:T-C0}\\
\|\hat T_\gamma'-T_\gamma'\|_{C^0(I_0)}
&\le C_\gamma^{(1)}\eps_F,\label{eq:T-C1}
\end{align}
where \(C_\gamma^{(1)}\) is explicit in \(C_0,C_1,L_0,L_1\) and \(\ell(\gamma)\). In particular, a bound \(\eps_T\le C_\gamma\eps_F\) is available on every fixed admissible loop.
\end{lemma}

\begin{proof}
Writing \(G=(A-yB)/(y-q)\), the identities
\[
\partial_yG=-\frac{A-qB}{(y-q)^2},
\qquad
\partial_y^2G=2\frac{A-qB}{(y-q)^3}
\]
give the first two estimates. Let $\eta_*:=\eta/2$ denote the common lower separation margin after applying Lemma~\ref{lem:separation-persistence}. For the exact--perturbed pair one may take the common Lipschitz constant
\begin{equation}
L_0=
\frac{1}{\eta_*^2}
\max\left\{
\sup_M\bigl[|u_F|(1+|q_F|\|S\|_{\mathrm{op}})\bigr],
\sup_M\bigl[|\hat u_F|(1+|\hat q_F|\|S\|_{\mathrm{op}})\bigr]
\right\}.
\label{eq:L0}
\end{equation}
Because \(\eps_F\le\eps_0\), this quantity is bounded uniformly in terms of the exact field bounds, \(\eps_0\), \(\eta\), and \(\|S\|_{\mathrm{op}}\). The perturbation estimates follow by differentiating the rational function \((A-yB)/(y-q)\) with respect to \(A,B,q\) on the common separation region.

If \(e=\hat y-y\) denotes the difference of lifts with the same initial value, Gronwall's inequality gives
\[
|e(1)|\le C_0\ell(\gamma)e^{L_0\ell(\gamma)}\eps_F.
\]
For \(Z=\partial_a y\) and \(\hat Z=\partial_a\hat y\),
\[
Z'=G_y(t,y)Z,
\qquad
\hat Z'=\hat G_y(t,\hat y)\hat Z.
\]
The exact variation satisfies \(|Z(t)|\le e^{L_0\ell(\gamma)}\). Subtracting the two equations and using
\[
|\hat G_y(t,\hat y)-G_y(t,y)|
\le L_1|\dot\gamma|\,|e|+C_1|\dot\gamma|\eps_F
\]
followed by a second Gronwall estimate proves \eqref{eq:T-C1}.
\end{proof}

\subsection{Global propagation in the non-neutral regime}

\begin{corollary}[Global non-neutral reconstruction stability]
\label{cor:globalnonneutral}
Assume the hypotheses of Theorem~\ref{thm:nonneutral} and Lemma~\ref{lem:reduction}, uniformly along paths from \(c_0\) to \(c\). Let \(\hat d\) be reconstructed by transport with \(\hat d(c_0)=\hat a\). Then
\[
|\hat d(c)-d(c)|
\le
e^{L_0\ell(\gamma_c)}
\left(|\hat a-a_*|+C_0\ell(\gamma_c)\eps_F\right),
\]
with \(L_0\) as in \eqref{eq:L0}. If \(M\) is connected and
\[
D_M:=\operatorname{diam}_{d_M}(M)
\]
is its intrinsic Riemannian diameter, then
\[
\|\hat d-d\|_{C^0(M)}
\le
e^{L_0D_M}
\left(
\frac{4\eps_T}{\delta_\gamma}+C_0D_M\eps_F
\right).
\]
\end{corollary}

\section{The holonomy variational jet}
\label{sec:variational}

\subsection{First variation}

Define
\begin{equation}
Z_1(t):=\partial_a y_a(t)\big|_{a=a_*}.
\end{equation}
Then
\begin{equation}
Z_1'(t)=G_y(t,y_*(t))Z_1(t),\qquad Z_1(0)=1,
\end{equation}
so
\begin{equation}
Z_1(t)=\exp\left(\int_0^tG_y(s,y_*(s))\,ds\right).
\end{equation}
Hence
\begin{equation}
T_\gamma'(a_*)=Z_1(1),\qquad H'(a_*)=Z_1(1)-1.
\end{equation}
The neutral condition is therefore
\begin{equation}
Z_1(1)=1.
\end{equation}

\subsection{Second variation: exact formula}

Define
\begin{equation}
W(t):=\partial_a^2y_a(t)\big|_{a=a_*}.
\end{equation}
Differentiating the variational equation gives
\begin{equation}
W'=G_y(t,y_*)W+G_{yy}(t,y_*)Z_1^2,\qquad W(0)=0.
\end{equation}
Since \(Z_1'=G_yZ_1\),
\begin{equation}
\left(\frac{W}{Z_1}\right)'=G_{yy}(t,y_*)Z_1.
\end{equation}
Therefore
\begin{equation}
W(1)=Z_1(1)\int_0^1G_{yy}(t,y_*(t))Z_1(t)\,dt.
\end{equation}

\begin{theorem}[Exact second holonomy variation]
\label{thm:secondvariation}
Assume that the horizontal transport is defined on a common fibre interval around \(a_*\) and is twice differentiable with respect to the initial anchor. Then
\begin{equation}
H''(a_*)=T_\gamma''(a_*)
=
T_\gamma'(a_*)\int_0^1G_{yy}(t,d(\gamma(t)))Z_1(t)\,dt.
\end{equation}
Equivalently,
\begin{equation}
H''(a_*)
=
2T_\gamma'(a_*)
\int_0^1
\frac{P(t)}{(d(\gamma(t))-q(t))^3}
Z_1(t)\,dt.
\end{equation}
In the neutral regime \(T_\gamma'(a_*)=1\), one has
\begin{equation}
\boxed{
H''(a_*)
=
2\int_0^1
\frac{P(t)}{(d(\gamma(t))-q(t))^3}
Z_1(t)\,dt.
}
\end{equation}
\end{theorem}

\begin{proof}
The formula follows from the variational equation for \(W\) and the identity
\[
(W/Z_1)'=G_{yy}Z_1.
\]
The explicit expression uses
\[
G_{yy}(t,y)=2\frac{P(t)}{(y-q(t))^3}.
\]
\end{proof}

\begin{remark}
Theorem \ref{thm:secondvariation} replaces the earlier assertion that second-order holonomy is ``controlled by \(d\eta_d\)''. The correct object is the vertical second variation of the nonlinear transport along the true horizontal section.
\end{remark}

\subsection{Higher-order holonomy jet}

For $n\ge1$, define
\[
Z_n(t):=\left.\partial_a^n y_a(t)\right|_{a=a_*}.
\]
The first two equations are
\[
Z_1'=G_yZ_1,
\qquad
Z_2'=G_yZ_2+G_{yy}Z_1^2.
\]
The general equation is most compactly written with the partial Bell polynomials $B_{n,r}$. We use the convention
\[
B_{n,r}(x_1,\ldots,x_{n-r+1})
=
\sum
\frac{n!}{j_1!\cdots j_{n-r+1}!}
\prod_{m=1}^{n-r+1}\left(\frac{x_m}{m!}\right)^{j_m},
\]
where the sum is over $j_m\ge0$ satisfying $\sum_m j_m=r$ and $\sum_m m j_m=n$.

\begin{theorem}[Finite holonomy variational hierarchy]
\label{thm:bell-jet}
Assume the horizontal transport is defined on a common fibre interval around $a_*$ and that $G_\gamma(t,\cdot)$ is $C^p$ there. Then, for $1\le n\le p$,
\begin{equation}
Z_n'
=
\sum_{r=1}^{n}
G_{y^r}(t,y_*(t))
B_{n,r}(Z_1,\ldots,Z_{n-r+1}),
\label{eq:bell-variation}
\end{equation}
with
\[
Z_1(0)=1,
\qquad
Z_n(0)=0\quad(n\ge2).
\]
At the endpoint,
\[
H'(a_*)=Z_1(1)-1,
\qquad
H^{(n)}(a_*)=Z_n(1)\quad(n\ge2).
\]
\end{theorem}

\begin{proof}
Differentiate the identity $y_a'=G_\gamma(t,y_a)$ $n$ times with respect to $a$. Fa\`a di Bruno's formula gives \eqref{eq:bell-variation}; the initial conditions follow from $y_a(0)=a$. Since $T_\gamma(a)=y_a(1)$ and $H=T_\gamma-\Id$, the endpoint identities follow.
\end{proof}

For the third derivative the recursion is explicit:
\[
Z_3'
=G_yZ_3+3G_{yy}Z_1Z_2+G_{yyy}Z_1^3,
\qquad
Z_3(0)=0.
\]
Dividing by $Z_1$ and using $Z_1'=G_yZ_1$ gives the following useful formula.

\begin{corollary}[Exact third holonomy variation]
\label{cor:thirdvariation}
Under the hypotheses of Theorem~\ref{thm:bell-jet},
\begin{equation}
H'''(a_*)
=
T_\gamma'(a_*)
\int_0^1
\left(
3G_{yy}(t,y_*)Z_2(t)
+G_{yyy}(t,y_*)Z_1(t)^2
\right)dt,
\label{eq:thirdvariation}
\end{equation}
where
\[
G_{yyy}(t,y)
=-6\frac{P(t)}{(y-q(t))^4}.
\]
\end{corollary}

The holonomy jet is therefore the endpoint jet of the vertical variational hierarchy:
\[
\J_\gamma(d;a_*)
=
\bigl(H'(a_*),H''(a_*),H'''(a_*),\ldots\bigr).
\]
The isolation order is the first nonzero component of this jet.

\subsection{Finite-dimensional holonomy unfoldings and Malgrange preparation}
\label{subsec:malgrange}

The finite-order condition is not only a Taylor-order statement.  For smooth finite-dimensional families of admissible return data it has an exact local preparation in the sense of Malgrange.  We use the standard smooth preparation theorem; see Golubitsky--Guillemin \cite[Chap.~IV]{GolubitskyGuillemin1973}.  The restriction to finite-dimensional parameter families is essential for this direct application: the full space of return-map perturbations is infinite-dimensional and is not covered by the classical theorem without additional functional-analytic structure.

Let $\Lambda\subset\mathbb R^m$ be a neighbourhood of $0$.  Consider a $C^\infty$ family of return components
\[
\lambda\longmapsto (u_\lambda,q_\lambda),\qquad \lambda\in\Lambda,
\]
along a fixed closed loop $\gamma$, and let $T_{\gamma,\lambda}$ denote the corresponding loop transport.  Assume that the horizontal lifts are defined on a common fibre neighbourhood of $a_*$ and remain uniformly separated from the singular set.  Standard smooth dependence of ODE solutions on parameters and initial data then gives a jointly smooth displacement
\[
H(\lambda,a):=T_{\gamma,\lambda}(a)-a.
\]

\begin{theorem}[Malgrange preparation of a finite-order holonomy singularity]
\label{thm:malgrange-holonomy}
Assume that, at $\lambda=0$,
\[
H(0,a_*)=0,
\qquad
\partial_a^jH(0,a_*)=0\quad(1\le j<k),
\qquad
\partial_a^kH(0,a_*)\ne0.
\]
Set $\xi=a-a_*$.  Then there exist neighbourhoods $\Lambda_0\ni0$ and $J_0\ni0$, a smooth function
\[
U\in C^\infty(\Lambda_0\times J_0),
\qquad U(0,0)\ne0,
\]
and smooth coefficient functions $c_0,\ldots,c_{k-1}\in C^\infty(\Lambda_0)$ with
\[
c_j(0)=0\qquad(0\le j<k),
\]
such that, after shrinking the neighbourhoods if necessary,
\begin{equation}
H(\lambda,a_*+\xi)
=
U(\lambda,\xi)
\left(
\xi^k+\sum_{j=0}^{k-1}c_j(\lambda)\xi^j
\right).
\label{eq:malgrange-holonomy}
\end{equation}
In particular, $U$ is nonvanishing on $\Lambda_0\times J_0$, and the local real anchors of the perturbed holonomy equation $H(\lambda,a)=0$ are exactly the real roots of the prepared polynomial
\[
P_\lambda(\xi)
:=
\xi^k+\sum_{j=0}^{k-1}c_j(\lambda)\xi^j.
\]
\end{theorem}

\begin{proof}
The uniform separation assumption removes the pole of the rational transport field from a common neighbourhood of the relevant solution graphs.  Since the family $(u_\lambda,q_\lambda)$ is $C^\infty$, standard parameter dependence for ODEs implies that $(\lambda,a)\mapsto T_{\gamma,\lambda}(a)$, and hence $H$, is $C^\infty$ jointly near $(0,a_*)$.  In the coordinate $\xi=a-a_*$, the displayed derivative conditions say precisely that the smooth germ $(\lambda,\xi)\mapsto H(\lambda,a_*+\xi)$ is regular of order $k$ in the distinguished variable $\xi$.  The Malgrange preparation theorem therefore yields a smooth unit times a monic polynomial of degree $k$ in $\xi$, with smooth parameter-dependent coefficients.  Because the unperturbed germ has multiplicity exactly $k$, the prepared coefficients satisfy $c_j(0)=0$.  Shrinking the neighbourhoods makes the unit nonvanishing throughout.  The last assertion follows immediately from \eqref{eq:malgrange-holonomy}.
\end{proof}

\begin{corollary}[Quadratic neutral unfolding and discriminant]
\label{cor:malgrange-k2}
If $k=2$, then locally
\[
H(\lambda,a_*+\xi)
=U(\lambda,\xi)\bigl(\xi^2+c_1(\lambda)\xi+c_0(\lambda)\bigr).
\]
Hence the local real-anchor structure is governed by the discriminant
\[
\Delta(\lambda)=c_1(\lambda)^2-4c_0(\lambda).
\]
For sufficiently small $\lambda$, $\Delta>0$ gives two nearby simple real anchors, $\Delta=0$ a multiple real anchor, and $\Delta<0$ no nearby real anchor.  This is the prepared-polynomial counterpart of the even-order root-splitting alternatives in Proposition~\ref{prop:parity}.
\end{corollary}

\begin{proposition}[First-order coefficients of a quadratic prepared unfolding]
\label{prop:malgrange-k2-coeff}
Assume the hypotheses of Corollary~\ref{cor:malgrange-k2}. Write all derivatives below at $(\lambda,a)=(0,a_*)$, and set
\[
A:=H_{aa}(0,a_*)\ne0,
\qquad
B:=H_{aaa}(0,a_*).
\]
Then, as $\lambda\to0$,
\begin{align}
c_0(\lambda)
&=\frac{2}{A}\,\lambda\cdot\nabla_\lambda H+O(|\lambda|^2),\label{eq:c0-malgrange}\\
c_1(\lambda)
&=\frac{2}{A}\,\lambda\cdot\nabla_\lambda H_a
-\frac{2B}{3A^2}\,\lambda\cdot\nabla_\lambda H
+O(|\lambda|^2),\label{eq:c1-malgrange}\\
\Delta(\lambda)
&=-\frac{8}{A}\,\lambda\cdot\nabla_\lambda H+O(|\lambda|^2).\label{eq:delta-malgrange}
\end{align}
In particular, to first order the discriminant is controlled by the parameter derivative of the holonomy residual at the singular anchor, whereas the linear coefficient $c_1$ also feels the cubic anchor jet through $B$ unless $\nabla_\lambda H(0,a_*)=0$.
\end{proposition}

\begin{proof}
Write $\xi=a-a_*$ and
\[
H(\lambda,a_*+\xi)=U(\lambda,\xi)\bigl(\xi^2+c_1(\lambda)\xi+c_0(\lambda)\bigr).
\]
At $\lambda=0$, comparison with the Taylor expansion of $H(0,a_*+\xi)$ gives
\[
U(0,0)=\frac{A}{2},
\qquad
U_\xi(0,0)=\frac{B}{6}.
\]
Setting $\xi=0$ and differentiating with respect to $\lambda$ yields
\[
\nabla_\lambda c_0(0)=\frac{2}{A}\nabla_\lambda H(0,a_*).
\]
Differentiating first in $\xi$ and then in $\lambda$, and using $c_0(0)=c_1(0)=0$, gives
\[
\nabla_\lambda H_a
=U_\xi(0,0)\nabla_\lambda c_0(0)+U(0,0)\nabla_\lambda c_1(0),
\]
which is equivalent to \eqref{eq:c1-malgrange}. Finally $c_1(\lambda)^2=O(|\lambda|^2)$, so
\[
\Delta(\lambda)=c_1(\lambda)^2-4c_0(\lambda)
=-\frac{8}{A}\lambda\cdot\nabla_\lambda H+O(|\lambda|^2).
\]
\end{proof}

\begin{remark}[Structural versus quantitative information]
The preparation theorem does not replace the explicit H\"older estimates of Section~\ref{sec:finiteorder}.  It provides an exact local structural reduction for smooth finite-dimensional unfoldings, while Theorem~\ref{thm:orderk} quantifies how far a perturbed root can move under norm-bounded perturbations and Theorem~\ref{thm:finitejet-stability} connects those perturbations to noisy return data.  Thus the two viewpoints are complementary:
\[
\text{holonomy jet}
\Longrightarrow
\text{prepared polynomial}
\Longrightarrow
\text{discriminant/root geometry},
\]
with the stability theory supplying the corresponding conditioning estimates.
\end{remark}

\subsection{Analytic dependence and the flat regime}

The rational fibre dependence of $G_\gamma$ gives more structure than $C^\infty$ regularity. Analytic dependence of flows on initial conditions for time-varying analytic vector fields is standard; see, for example, \cite{JafarpourLewis2014}. We include the short argument in the present scalar setting.

\begin{theorem}[Analyticity of the holonomy map]
\label{thm:analytic-holonomy}
Let $J\Subset\mathbb R$ be an interval of initial anchors such that, for every $a\in J$, the horizontal lift $y_a(t)$ exists for $0\le t\le1$ and
\[
y_a(t)-q(t)\ge\eta>0
\]
uniformly in $(a,t)\in J\times[0,1]$. Then $a\mapsto T_\gamma(a)$ is real analytic on $J$.
Consequently, if $a_*\in J$, $H(a_*)=0$, and
\[
H^{(n)}(a_*)=0\qquad(n\ge1),
\]
then there exists a neighbourhood $U\subset J$ of $a_*$ on which
\[
T_\gamma(a)=a.
\]
Thus infinite isolation order for one loop is equivalent locally to identity holonomy for that loop.
\end{theorem}

\begin{proof}
Fix $a_0\in J$. By the uniform separation and continuous dependence of solutions, the exact solution graph through $a_0$ is contained in a compact real strip that has positive distance from the pole $y=q(t)$. The rational function $G_\gamma(t,z)$ therefore extends holomorphically in $z$ to a complex neighbourhood of that strip, uniformly for $t\in[0,1]$. On a sufficiently short time interval, Picard iteration for the complexified ODE maps a common complex neighbourhood of $a_0$ into itself; every iterate is holomorphic in the initial value and the iterates converge uniformly. Hence the local flow is holomorphic in the initial value. A finite composition of such short-time flow maps covers $[0,1]$, so $T_\gamma$ is holomorphic near $a_0$ and therefore real analytic on $J$.

Since $H=T_\gamma-\Id$ is real analytic, vanishing of $H(a_*)$ and all derivatives at $a_*$ forces its Taylor series to vanish identically in a neighbourhood of $a_*$. The identity theorem gives $H\equiv0$ there.
\end{proof}

\section{Finite-order anchor isolation and parity-sensitive existence}
\label{sec:finiteorder}

The scalar exponent \(1/k\) follows from the Taylor structure of a zero of multiplicity \(k\). In the present inverse problem, the substantive point is that the multiplicity is the first nonzero component of the holonomy jet derived in Section~\ref{sec:variational}.

\begin{theorem}[Order-\(k\) isolation]
\label{thm:orderk}
Let \(k\ge2\), let \(H\in C^{k+1}(\bar I)\), and assume
\[
H(a_*)=0,\qquad H^{(j)}(a_*)=0\quad(1\le j<k),
\qquad
m_k:=|H^{(k)}(a_*)|>0.
\]
Set
\[
M_{k+1}:=\|H^{(k+1)}\|_{C^0(\bar I)}.
\]
If \(M_{k+1}>0\), set
\[
r_k:=\frac{m_k}{4M_{k+1}}
\]
and assume \(\overline{B_{r_k}(a_*)}\subset I\). Let \(\hat H=H+E\), with \(\|E\|_{C^0(\bar I)}\le\eps\), and suppose \(\hat H(\hat a)=0\). Assume either
\begin{enumerate}[label=(\alph*)]
\item \(M_{k+1}>0\) and \(|\hat a-a_*|\le r_k\), or
\item \(M_{k+1}=0\) and \(\hat a\in I\).
\end{enumerate}
Then
\begin{equation}
|\hat a-a_*|
\le
\kappa_k^{1/k}
\left(\frac{\eps}{m_k}\right)^{1/k},
\qquad
\kappa_k:=\frac{4(k+1)!}{4k+3}.
\label{eq:orderk}
\end{equation}
\end{theorem}

\begin{proof}
Taylor's formula gives
\[
H(\hat a)
=\frac{H^{(k)}(a_*)}{k!}(\hat a-a_*)^k+R_{k+1}.
\]
If \(M_{k+1}>0\), the localization hypothesis implies
\[
|R_{k+1}|
\le\frac{M_{k+1}}{(k+1)!}|\hat a-a_*|^{k+1}
\le\frac{m_k}{4(k+1)!}|\hat a-a_*|^k.
\]
Since \(\hat H(\hat a)=0\), one has \(|H(\hat a)|\le\eps\), and therefore
\[
\frac{m_k(4k+3)}{4(k+1)!}|\hat a-a_*|^k\le\eps.
\]
This is \eqref{eq:orderk}. If \(M_{k+1}=0\), the Taylor remainder vanishes identically on \(I\); the same bound remains valid (and is slightly weaker than the exact remainder-free estimate).
\end{proof}

\begin{proposition}[Parity-sensitive local existence]
\label{prop:parity}
Let $k\ge2$, let $H\in C^{k+1}(\bar I)$, and assume
\[
H(a_*)=0,
\qquad
H^{(j)}(a_*)=0\quad(1\le j<k),
\qquad
m_k:=|H^{(k)}(a_*)|>0.
\]
Set $M_{k+1}:=\|H^{(k+1)}\|_{C^0(\bar I)}$ and let $\hat H$ satisfy
\[
\|\hat H-H\|_{C^0(I)}\le\eps.
\]
Choose $\rho>0$ such that
\[
[a_*-\rho,a_*+\rho]\subset I
\]
and, when $M_{k+1}>0$,
\[
\rho\le\frac{(k+1)m_k}{2M_{k+1}}.
\]
If
\[
\eps<\frac{m_k}{4k!}\rho^k,
\]
then:
\begin{enumerate}[label=(\roman*)]
\item if $k$ is odd, $\hat H$ has at least one zero in $(a_*-\rho,a_*+\rho)$;
\item if $k$ is even and $\hat H(a_*)H^{(k)}(a_*)<0$, then $\hat H$ has at least one zero in each of $(a_*-\rho,a_*)$ and $(a_*,a_*+\rho)$;
\item if $k$ is even and $\hat H(a_*)=0$, then $a_*$ itself remains a zero.
\end{enumerate}
\end{proposition}

\begin{proof}
Taylor's formula and the choice of \(\rho\) give
\[
|H(a_*\pm\rho)|\ge\frac{m_k}{2k!}\rho^k.
\]
The perturbation bound preserves the endpoint signs. For odd \(k\), these signs are opposite, so the intermediate value theorem gives a zero. For even \(k\), the endpoint signs agree with \(H^{(k)}(a_*)\); the stated sign condition at \(a_*\) gives one zero on each side. The last case is immediate.
\end{proof}

\begin{corollary}[Finite-order stability hierarchy]
\label{cor:hierarchy}
Let \(k_\gamma(a_*)<\infty\).
\begin{enumerate}[label=(\roman*)]
\item If \(k_\gamma(a_*)=1\), Theorem~\ref{thm:nonneutral} gives first-order stability \(|\hat a-a_*|=O(\eps_T)\).
\item If \(k_\gamma(a_*)=k\ge2\), every perturbed zero satisfying the localization hypothesis of Theorem~\ref{thm:orderk} obeys
\[
|\hat a-a_*|
\le\kappa_k^{1/k}
\left(\frac{\eps}{m_k}\right)^{1/k}.
\]
Existence is automatic for sufficiently small perturbations when \(k\) is odd, and for even \(k\) follows under the sign alternatives of Proposition~\ref{prop:parity}.
\end{enumerate}
If \(k_\gamma(a_*)=+\infty\), no finite-order location estimate follows from that loop alone.
\end{corollary}

\begin{remark}[Size of the quantitative constants]
For the pure local model
\[
H(a)=\frac{m_k}{k!}(a-a_*)^k,
\qquad
\hat H=H-\varepsilon,
\]
the exact root displacement is
\[
|\hat a-a_*|=(k!)^{1/k}\left(\frac{\varepsilon}{m_k}\right)^{1/k}.
\]
Theorem~\ref{thm:orderk} gives instead the coefficient $\kappa_k^{1/k}$. For $k=2,3,4$ the exact pure-root coefficients are approximately $1.4142,1.8171,2.2134$, while the theorem gives $1.4771,1.8566,2.2419$. Thus the explicit Taylor-remainder constants are quantitative rather than optimal, but are already close to the pure-root asymptotic constants at low order.
\end{remark}

\section{Stability and observable finite-order coefficients}
\label{sec:observable-k}

The previous section assumes control of the holonomy displacement. We now close the inverse-problem chain by showing that perturbations of the observable return components control every fixed finite vertical holonomy jet.

\begin{theorem}[Finite-jet stability under return-data perturbations]
\label{thm:finitejet-stability}
Fix an integer $p\ge0$, a closed $C^1$ loop $\gamma$, a compact anchor interval $I_0$, and a perturbation ceiling $\varepsilon_0>0$. Let $(u_F,q_F)$ and $(\hat u_F,\hat q_F)$ be exact and perturbed return components satisfying
\[
\|\hat u_F-u_F\|_{C^0(\gamma)}+\|\hat q_F-q_F\|_{C^0(\gamma)}\le\varepsilon_F\le\varepsilon_0.
\]
Assume that the exact horizontal lifts are uniformly separated on $I_0$ with margin $\eta>0$ in the sense of Definition~\ref{def:uniform-separation}. Shrinking $\varepsilon_0$ if necessary, Lemma~\ref{lem:separation-persistence} supplies existence, a common bounded fibre region, and perturbed separation with margin $\eta/2$. Then, for every $0\le j\le p$, there is a constant $C_{p,\gamma,j}$ depending only on $p$, $\eta$, $\varepsilon_0$, the exact geometric bounds, $I_0$ and $\gamma$ such that
\begin{equation}
\|\hat T_\gamma^{(j)}-T_\gamma^{(j)}\|_{C^0(I_0)}
\le
C_{p,\gamma,j}\varepsilon_F.
\label{eq:finitejet-stability}
\end{equation}
In particular,
\[
\|\hat H^{(j)}-H^{(j)}\|_{C^0(I_0)}
\le C_{p,\gamma,j}\varepsilon_F
\qquad(j\ge1),
\]
while the same conclusion for $j=0$ follows from $\hat H-H=\hat T_\gamma-T_\gamma$.
\end{theorem}

\begin{proof}
For $r\ge1$, the exact fibre derivatives are
\[
\partial_y^rG_\gamma(t,y)=(-1)^r r!\frac{P(t)}{(y-q(t))^{r+1}},
\]
and the perturbed derivatives have the same form with hatted data. Uniform separation and boundedness imply, for $0\le r\le p+1$, uniform bounds
\[
|\partial_y^rG_\gamma|+|\partial_y^r\hat G_\gamma|\le C_r
\]
and, by direct subtraction of the rational expressions,
\[
|\partial_y^r\hat G_\gamma-\partial_y^rG_\gamma|
\le C_r'\varepsilon_F.
\]
The derivative at order $p+1$ is used when comparing $\partial_y^p\hat G_\gamma(t,\hat y)$ and $\partial_y^p\hat G_\gamma(t,y)$ by the mean-value estimate in the fibre variable. The case $j=0$ is the standard Gronwall comparison of the two transport equations. For $j\ge1$, let $Z_j=\partial_a^jy_a$ and $\hat Z_j=\partial_a^j\hat y_a$. Theorem~\ref{thm:bell-jet} gives triangular Bell-polynomial systems for $Z_j$ and $\hat Z_j$. Uniform bounds for the lower-order variations follow inductively from Gronwall. Subtracting the two equations at order $j$ yields a linear term in $\hat Z_j-Z_j$ plus a forcing bounded by $C\varepsilon_F$ and lower-order differences, already controlled by the induction hypothesis. A further Gronwall estimate gives
\[
\|\hat Z_j-Z_j\|_{C^0(I_0\times[0,1])}\le C_{p,\gamma,j}\varepsilon_F.
\]
Evaluating at $t=1$ proves \eqref{eq:finitejet-stability}.
\end{proof}

\begin{remark}
Theorem~\ref{thm:finitejet-stability} is specific to the vertical anchor jet. It requires only $C^0$ perturbation of the observable fields along the loop because all anchor derivatives are derivatives with respect to the scalar fibre variable; no spatial differentiation of the noisy return data is introduced.
\end{remark}

The theoretical coefficient $m_k=|H^{(k)}(a_*)|$ is evaluated at the unknown true anchor. It is therefore useful to quantify how accurately it can be replaced by a coefficient computed at the reconstructed anchor.

\begin{proposition}[Observable holonomy coefficient]
\label{prop:observable-mk}
Let \(k\ge2\). Assume \(H\in C^{k+1}(\bar I)\), \(\hat H\in C^k(\bar I)\), and
\[
\|\hat H^{(k)}-H^{(k)}\|_{C^0(\bar I)}\le\eps_k.
\]
Set
\[
\hat m_k:=|\hat H^{(k)}(\hat a)|,
\qquad
m_k:=|H^{(k)}(a_*)|,
\qquad
M_{k+1}:=\|H^{(k+1)}\|_{C^0(\bar I)}.
\]
Then
\begin{equation}
|\hat m_k-m_k|
\le\eps_k+M_{k+1}|\hat a-a_*|.
\label{eq:mk-proximity}
\end{equation}
If the right-hand side is at most \(\theta m_k\) for some \(0<\theta<1\), then
\[
\frac{1}{1+\theta}\hat m_k\le m_k\le\frac{1}{1-\theta}\hat m_k.
\]
Consequently, whenever Theorem~\ref{thm:orderk} applies,
\begin{equation}
|\hat a-a_*|
\le
\bigl((1+\theta)\kappa_k\bigr)^{1/k}
\left(\frac{\eps}{\hat m_k}\right)^{1/k}.
\label{eq:observable-mk}
\end{equation}
\end{proposition}

\begin{proof}
The reverse triangle inequality gives
\[
|\hat m_k-m_k|
\le
|\hat H^{(k)}(\hat a)-H^{(k)}(\hat a)|
+|H^{(k)}(\hat a)-H^{(k)}(a_*)|,
\]
which is \eqref{eq:mk-proximity}. The two-sided comparison of \(m_k\) and \(\hat m_k\) is immediate, and substituting \(m_k\ge\hat m_k/(1+\theta)\) into Theorem~\ref{thm:orderk} gives \eqref{eq:observable-mk}.
\end{proof}

For \(k=2\) and \(\theta=1/2\), the observable estimate becomes
\begin{equation}
|\hat a-a_*|
\le
\left(\frac{36\eps}{11\hat m_2}\right)^{1/2},
\label{eq:observable-m2}
\end{equation}
whereas the theorem constant expressed with the unknown exact coefficient \(m_2\) is \(\sqrt{24/11}\). The leading scalar asymptotic coefficient for a pure quadratic root is \(\sqrt{2}\), so the displayed theorem constant is quantitative rather than optimal.

\begin{corollary}[Observable coefficients from return data]
\label{cor:observable-from-F}
Under the hypotheses of Theorem~\ref{thm:finitejet-stability}, for every fixed $k\ge2$,
\[
|\hat m_k-m_k|
\le
C_{k,\gamma}\varepsilon_F
+M_{k+1}|\hat a-a_*|.
\]
Hence whenever the right-hand side is at most $\theta m_k$ with $0<\theta<1$, the observable estimate \eqref{eq:observable-mk} applies with a coefficient computable from the perturbed holonomy jet.
\end{corollary}

\section{Multi-loop finite-order aggregation}
\label{sec:multiloop}

A single loop may have a weak or high-order signal. Several loops sharing the same true anchor can be combined deterministically without assuming statistical independence.

Let $\gamma_1,\ldots,\gamma_r$ be closed loops based at the same $c_0$, with holonomy displacements $H_j=T_{\gamma_j}-\Id$. Set
\[
k_*:=\min_j k_{\gamma_j}(a_*).
\]
Assume $k_*<\infty$ and set
\[
J_*:=\{j:k_{\gamma_j}(a_*)=k_*\}.
\]
For $j\in J_*$ define
\[
v_j:=H_j^{(k_*)}(a_*),
\qquad
m_*:=\left(\sum_{j\in J_*}|v_j|^2\right)^{1/2}>0,
\]
and
\[
M_*:=\left(\sum_{j\in J_*}\|H_j^{(k_*+1)}\|_{C^0(\bar I)}^2\right)^{1/2}.
\]

\begin{theorem}[Deterministic multi-loop finite-order bound]
\label{thm:multiloop}
Assume $H_j(a_*)=0$ and $H_j^{(\ell)}(a_*)=0$ for $1\le\ell<k_*$ and every $j\in J_*$. Let
\[
\mathbf H=(H_j)_{j\in J_*},
\qquad
\hat{\mathbf H}=(\hat H_j)_{j\in J_*},
\]
and suppose
\[
\sup_{a\in I}\|\hat{\mathbf H}(a)-\mathbf H(a)\|_2\le\varepsilon_*.
\]
Let $\hat a\in I$ be an approximate common anchor satisfying
\[
\|\hat{\mathbf H}(\hat a)\|_2\le\tau.
\]
If $M_*>0$, set $r_*=m_*/(4M_*)$, assume $\overline{B_{r_*}(a_*)}\subset I$, and suppose $|\hat a-a_*|\le r_*$. If $M_*=0$, assume only $\hat a\in I$. Then
\begin{equation}
|\hat a-a_*|
\le
\kappa_{k_*}^{1/k_*}
\left(\frac{\varepsilon_*+\tau}{m_*}\right)^{1/k_*},
\qquad
\kappa_k=\frac{4(k+1)!}{4k+3}.
\label{eq:multiloop}
\end{equation}
\end{theorem}

\begin{proof}
Write $s=\hat a-a_*$. Vector Taylor expansion gives
\[
\mathbf H(a_*+s)
=
\frac{s^{k_*}}{k_*!}\,v+\mathbf R_{k_*+1}(s),
\qquad
v=(v_j)_{j\in J_*},
\]
with
\[
\|\mathbf R_{k_*+1}(s)\|_2
\le
\frac{M_*}{(k_*+1)!}|s|^{k_*+1}.
\]
Under the localization hypothesis,
\[
\|\mathbf H(a_*+s)\|_2
\ge
\frac{m_*(4k_*+3)}{4(k_*+1)!}|s|^{k_*}.
\]
On the other hand,
\[
\|\mathbf H(\hat a)\|_2
\le
\|\mathbf H(\hat a)-\hat{\mathbf H}(\hat a)\|_2
+\|\hat{\mathbf H}(\hat a)\|_2
\le\varepsilon_*+\tau.
\]
Combining the two estimates yields \eqref{eq:multiloop}. If $M_*=0$, the remainder vanishes and the same bound follows directly.
\end{proof}

\begin{remark}[Loop-design consequence]
The theorem gives a deterministic ordering principle: first minimize the isolation order $k_\gamma$; among loops of the minimal order, increase the collective coefficient $m_*$ while controlling loop length and transport conditioning. The residual $\tau$ allows numerically computed loops that do not possess an exactly common perturbed fixed point.
\end{remark}

\section{Global reconstruction stability at finite holonomy order}
\label{sec:globalneutral}

The local anchor estimate propagates through the same scalar transport equation used for exact reconstruction.

\begin{theorem}[Global finite-order reconstruction stability]
\label{thm:globalneutral}
Let \(k=k_\gamma(a_*)\ge2\) be finite. Assume:
\begin{enumerate}[label=(\roman*)]
\item the exact and perturbed transports satisfy the uniform separation and boundedness hypotheses of Lemma~\ref{lem:reduction} along paths connecting \(c_0\) to every point of \(M\);
\item \(\hat H\) has a zero \(\hat a\) satisfying the localization hypotheses of Theorem~\ref{thm:orderk} (for example, under the odd/even sufficient conditions of Proposition~\ref{prop:parity});
\item the holonomy perturbation obeys \(\eps\le C_\gamma\eps_F\) for the chosen identifying loop.
\end{enumerate}
Let
\[
D_M:=\operatorname{diam}_{d_M}(M).
\]
Then
\begin{equation}
\|\hat d-d\|_{C^0(M)}
\le
C_{\mathrm{glob}}
\left(\frac{\eps_F}{m_k}\right)^{1/k}
+C_{\mathrm{prop}}\eps_F,
\label{eq:global-k}
\end{equation}
where one may take
\[
C_{\mathrm{glob}}
=e^{L_0D_M}\kappa_k^{1/k}C_\gamma^{1/k},
\qquad
C_{\mathrm{prop}}=C_0D_Me^{L_0D_M},
\]
with \(L_0\) from \eqref{eq:L0}. Hence, for fixed \(m_k>0\),
\[
\|\hat d-d\|_{C^0(M)}=O(\eps_F^{1/k}).
\]
\end{theorem}

\begin{proof}
For a path \(\gamma_c\) from \(c_0\) to \(c\), the standard Gronwall comparison gives
\[
|\hat d(c)-d(c)|
\le
e^{L_0\ell(\gamma_c)}
\left(|\hat a-a_*|+C_0\ell(\gamma_c)\eps_F\right).
\]
Choose paths with \(\ell(\gamma_c)\le D_M\). By Theorem~\ref{thm:orderk},
\[
|\hat a-a_*|
\le\kappa_k^{1/k}
\left(\frac{\eps}{m_k}\right)^{1/k}
\le\kappa_k^{1/k}C_\gamma^{1/k}
\left(\frac{\eps_F}{m_k}\right)^{1/k}.
\]
Substitution yields \eqref{eq:global-k}.
\end{proof}

\section{Geometric interpretation and remaining geometric questions}
\label{sec:geometry}

The exact second- and third-variation formulas give rigorous geometric signals beyond the linearized multiplier: $H''(a_*)$ and $H'''(a_*)$ are integrals of vertical derivatives of the nonlinear return transport along the true horizontal section. Theorem~\ref{thm:bell-jet} extends this to every finite order, while Theorem~\ref{thm:analytic-holonomy} settles the local analytic meaning of the flat regime for a fixed loop. The material below is deliberately stated as geometric perspective and open problems; no unproved small-loop curvature expansion or genericity assertion is used in the theorem chain.

\subsection{Connection-curvature perspective}

The inverse transport can be written as the horizontal distribution of the one-form
\[
\Theta_F=dz-\omega_F(c,z)
\]
on the scalar bundle $E_F=\{(c,z):z>q_F(c)\}$, as in the nonlinear-holonomy formulation of \cite{BarkatouElMorsalaniInv2026}; the connection language is standard, see \cite{KobayashiNomizu1963}. If
\[
X_i=\partial_i+\omega_i(c,z)\partial_z,
\]
then
\[
[X_i,X_j]_{\mathrm{vert}}
=
\bigl(\partial_i\omega_j-\partial_j\omega_i
+\omega_i\partial_z\omega_j
-\omega_j\partial_z\omega_i\bigr)\partial_z.
\]
This suggests the scalar curvature two-form
\[
\mathcal K_F=d_M\omega_F+\omega_F\wedge\partial_z\omega_F.
\]
A remaining invariant problem is to express small-loop holonomy-jet coefficients in terms of the vertical curvature jet $\partial_z^r\mathcal K_F$. Such a formula would connect the finite-order scalar singularity theory proved here to intrinsic local geometry of the inverse connection.

\subsection{Spherical and non-umbilic perturbations}

For a perturbation $d_\lambda=d_0+\lambda\phi$, the true anchor is
\[
a_*(\lambda)=d_\lambda(c_0)=d_0+\lambda\phi(c_0),
\]
and the true horizontal section is $y_*(t;\lambda)=d_\lambda(\gamma(t))$. Expansions in $\lambda$ must therefore be performed along this moving section.

Likewise, a compact convex surface in $\mathbb R^3$ should not be modelled globally by two distinct constant principal curvatures. For a genuine non-umbilic core, such as an ellipsoid, an invariant leading coefficient must include the Levi-Civita connection terms of the principal frame. Deriving such a coefficient and proving generic non-vanishing remain open geometric problems. These questions are distinct from the vertical anchor-jet hierarchy, which is exact once the observed transport is fixed.

\subsection{Neutral tuning}

If a one-parameter family satisfies
\[
H'(a_*;\lambda)=D_1\lambda+O(\lambda^2),
\qquad D_1\ne0,
\]
then the sign is fixed for sufficiently small positive $\lambda$; a nonzero linear coefficient prevents nearby neutral tuning. A generic search for neutral configurations should therefore use either a two-parameter family $d_{\lambda,\mu}=d_0+\lambda\phi+\mu\psi$ or a specially designed one-parameter family. Section~\ref{sec:numerics} shows that exact neutral tuning can nevertheless be achieved explicitly in a finite Fourier family on the circular core.

\section{Numerical experiments}
\label{sec:numerics}

The numerical section has two purposes. First, it calibrates the finite-order scalar exponent without geometric complications. Second, it tests the actual return-map transport on genuine admissible geometries. All computations reported below were performed in double precision with Python~3.13.5, NumPy~2.3.5 and SciPy~1.17.0. The transport equations and their first two variational equations were integrated by DOP853 with relative tolerance $2\times10^{-11}$ (and $2\times10^{-12}$ for the exact baseline derivatives), absolute tolerance at most $10^{-13}$, and maximum step $0.02$. Roots were refined by Brent's method to relative and absolute tolerances of order $10^{-13}$. Geometric margins were first sampled on a grid of $20{,}000$ angles and then enclosed independently by outward-rounded interval arithmetic (\texttt{mpmath.iv}) on $2{,}000$ angular subintervals. A reference implementation accompanies the manuscript.

\subsection{Algebraic finite-order calibration}

For $H(a)=a^k$ at $a_*=0$ and the perturbation $\hat H_\varepsilon(a)=a^k-\varepsilon$, the relevant real root has exact magnitude $\varepsilon^{1/k}$ (for even $k$ we take the positive branch). Thus the exponents $1/2$, $1/3$ and $1/4$ are exact for $k=2,3,4$.

\begin{table}[ht]
\centering
\caption{Exact scalar calibration of the finite-order exponent.}
\begin{tabular}{@{}cccc@{}}
\toprule
order $k$ & exact displacement & predicted exponent & observed log--log slope\\
\midrule
$2$ & $\varepsilon^{1/2}$ & $1/2$ & $0.500000$\\
$3$ & $\varepsilon^{1/3}$ & $1/3$ & $0.333333$\\
$4$ & $\varepsilon^{1/4}$ & $1/4$ & $0.250000$\\
\bottomrule
\end{tabular}
\end{table}

This benchmark only tests scalar singular-root geometry. The next experiments use return data generated by actual normal-incidence geometry.

\subsection{Circular-core return geometry}

Let $C$ be the unit disk and write
\[
c(\theta)=(\cos\theta,\sin\theta),
\qquad
\tau(\theta)=(-\sin\theta,\cos\theta),
\qquad S=-I.
\]
For a positive thickness $d(\theta)$ set
\[
w(\theta)=\frac{d'(\theta)}{1+d(\theta)}.
\]
The outer point is $\Phi(\theta)=(1+d(\theta))c(\theta)$ and its inward normal is
\[
n(\theta)=\frac{-c(\theta)+w(\theta)\tau(\theta)}{\sqrt{1+w(\theta)^2}}.
\]
Writing the return decomposition as
\[
F(\theta)-c(\theta)=u(\theta)\tau(\theta)+q(\theta)c(\theta),
\]
one obtains
\begin{align}
D(\theta)&=1-w(\theta)^2\bigl((1+d(\theta))^2-1\bigr),\\
\sigma(\theta)&=\frac{1+d(\theta)-\sqrt{D(\theta)}}{1+w(\theta)^2},\\
q(\theta)&=d(\theta)-\sigma(\theta),
\qquad
u_F(\theta):=u(\theta)=\sigma(\theta)w(\theta).
\end{align}
The inverse transport around the positively oriented circle is therefore
\begin{equation}
y'(\theta)=\frac{u_F(\theta)(1+y(\theta))}{y(\theta)-q(\theta)},
\qquad 0\le\theta\le2\pi.
\label{eq:circle-transport}
\end{equation}
For the perturbation experiments we keep $u_F$ fixed and add smooth measurement error to the observed normal component,
\begin{equation}
\hat q_\rho(\theta)=q(\theta)-\rho\cos4\theta,
\qquad
\varepsilon_F=\|\hat q_\rho-q\|_\infty=\rho.
\label{eq:q-noise}
\end{equation}
The perturbed pair need not itself be an exact return map; it represents smooth observational error in the return components, which is precisely the setting of the stability theory.

\subsection{A genuine non-neutral example}

Take
\begin{equation}
d_{\mathrm{nn}}(\theta)
=0.8+0.2(\cos\theta+\sin\theta+\cos2\theta+\sin2\theta).
\label{eq:d-nonneutral}
\end{equation}
The true anchor at $\theta=0$ is $a_*=1.2$. Direct integration of \eqref{eq:circle-transport} and its variational equations gives
\[
H(a_*)=-5.6\times10^{-15},
\qquad
H'(a_*)=2.94769678\times10^{-2},
\qquad
H''(a_*)=-1.39481188\times10^{-1}.
\]
The geometry is well separated from degeneracy. A $20{,}000$-point grid gives
\[
\min d_{\mathrm{nn}}=0.38140,
\quad
\min D=0.50373,
\quad
\min\sigma=0.38140,
\quad
\max\langle n(\Phi),\nu(F)\rangle=-0.65104.
\]
The independent interval enclosure on $2{,}000$ angular cells yields the outward-rounded bounds
\[
d_{\mathrm{nn}}>0.38022,
\qquad
D>0.49950,
\qquad
\sigma>0.38022,
\qquad
\langle n(\Phi),\nu(F)\rangle<-0.64323.
\]
Thus the normal graph stays embedded and the first-hit discriminant is uniformly positive. The transversality is uniformly negative throughout the circle.

The perturbation response at the true anchor satisfies numerically
\[
\left.\partial_\rho\hat H_\rho(a_*)\right|_{\rho=0}=0.13825171.
\]
First-order root perturbation therefore predicts
\[
|\hat a-a_*|\sim
\frac{0.13825171}{0.0294769678}\,\rho
=4.69016\,\rho.
\]
At $\rho=10^{-8}$ the observed coefficient is $4.69017$, agreeing with the predicted leading constant to the displayed precision. For the perturbation \eqref{eq:q-noise}, Table~\ref{tab:nonneutral-num} shows the reconstructed anchor and global thickness errors. A least-squares fit over all nine noise levels $10^{-8}\le\rho\le10^{-4}$ gives slopes $1.00018$ and $1.00023$, respectively.

\begin{table}[ht]
\centering
\caption{Non-neutral geometric experiment.}
\label{tab:nonneutral-num}
\begin{tabular}{@{}ccc@{}}
\toprule
$\rho=\varepsilon_F$ & $|\hat a-a_*|$ & $\|\hat d-d\|_\infty$\\
\midrule
$10^{-8}$ & $4.69017\times10^{-8}$ & $9.02152\times10^{-8}$\\
$10^{-7}$ & $4.69017\times10^{-7}$ & $9.02153\times10^{-7}$\\
$10^{-6}$ & $4.69027\times10^{-6}$ & $9.02178\times10^{-6}$\\
$10^{-5}$ & $4.69130\times10^{-5}$ & $9.02430\times10^{-5}$\\
$10^{-4}$ & $4.70162\times10^{-4}$ & $9.04963\times10^{-4}$\\
\bottomrule
\end{tabular}
\end{table}

\begin{figure}[ht]
\centering
\begin{tikzpicture}
\begin{loglogaxis}[
width=.72\textwidth,height=.43\textwidth,
xlabel={$\varepsilon_F$},ylabel={error},
legend pos=north west,grid=major]
\addplot[mark=*] coordinates {
(1e-8,4.69016956e-8) (3e-8,1.40705000e-7) (1e-7,4.69017248e-7)
(3e-7,1.40705850e-6) (1e-6,4.69027484e-6) (3e-6,1.40715092e-5)
(1e-5,4.69130208e-5) (3e-5,1.40807626e-4) (1e-4,4.70161640e-4)};
\addlegendentry{$|\hat a-a_*|$}
\addplot[mark=square*] coordinates {
(1e-8,9.02151946e-8) (3e-8,2.70645451e-7) (1e-7,9.02153043e-7)
(3e-7,2.70647576e-6) (1e-6,9.02178182e-6) (3e-6,2.70670257e-5)
(1e-5,9.02430287e-5) (3e-5,2.70897385e-4) (1e-4,9.04963187e-4)};
\addlegendentry{$\|\hat d-d\|_\infty$}
\addplot[dashed,domain=1e-8:1e-4,samples=2] {6*x};
\addlegendentry{reference slope $1$}
\end{loglogaxis}
\end{tikzpicture}
\caption{Non-neutral return geometry: both local anchor error and global thickness error are first order in the return-data perturbation.}
\end{figure}
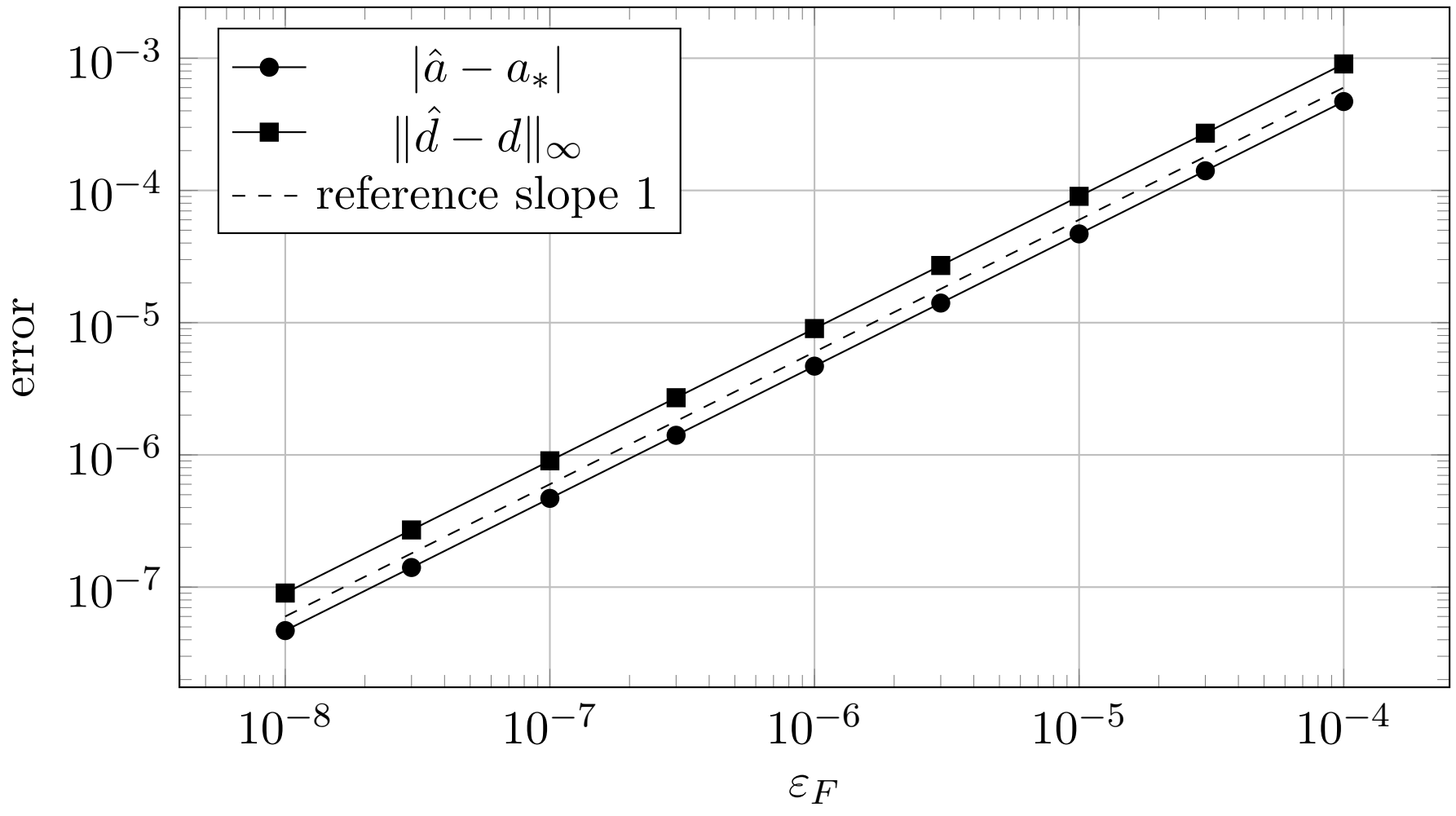

\subsection{A numerically tuned neutral example}

Consider the one-parameter Fourier family
\begin{align}
d_{\alpha}(\theta)=\;&0.8
+\alpha\cos\theta
+0.12866554\sin\theta\\
&-0.16790919\cos2\theta
+0.08267596\sin2\theta
-0.11676398\cos3\theta.
\end{align}
We tune the first coefficient by solving the scalar neutrality equation
\[
H_\alpha'\bigl(a_*(\alpha)\bigr)=0
\]
with Brent's method on $[0.14476636,0.14476637]$. This gives
\[
\alpha_*=0.14476636480382166.
\]
We then set
\begin{equation}
d_{\mathrm{neu}}:=d_{\alpha_*}.
\label{eq:d-neutral}
\end{equation}
At the base point $c_0=(1,0)$,
\[
a_*=d_{\mathrm{neu}}(0)=0.6600931948038217.
\]
The tuned loop transport is neutral to machine precision,
\[
H(a_*)=-2.2\times10^{-16},
\qquad
H'(a_*)=-8.9\times10^{-16},
\]
while the exact second-variation formula gives
\begin{equation}
H''(a_*)=3.13933613\times10^{-3}\ne0.
\label{eq:neutral-H2-num}
\end{equation}
The geometry again has substantial admissibility margins. The $20{,}000$-point grid gives
\[
\min d_{\mathrm{neu}}=0.58204,
\quad
\min D=0.56673,
\quad
\min\sigma=0.58204,
\quad
\max\langle n(\Phi),\nu(F)\rangle=-0.70167.
\]
The independent interval enclosure on $2{,}000$ angular cells gives
\[
d_{\mathrm{neu}}>0.58064,
\qquad
D>0.56375,
\qquad
\sigma>0.58064,
\qquad
\langle n(\Phi),\nu(F)\rangle<-0.69535.
\]
Hence the tuned normal graph is uniformly admissible for the numerical return construction. This is a numerically tuned second-order neutral holonomy generated by the return geometry, not an imposed scalar model.

For the perturbation \eqref{eq:q-noise}, the residual at the true anchor has the linear response
\[
\left.\partial_\rho\hat H_\rho(a_*)\right|_{\rho=0}
=-0.16119036.
\]
Combining this with \eqref{eq:neutral-H2-num} gives the quadratic normal-form prediction
\[
|\hat a-a_*|\sim
\sqrt{\frac{2(0.16119036)}{0.00313933613}}\,\sqrt{\rho}
=10.13364\,\sqrt{\rho}.
\]
At $\rho=10^{-8}$ the observed upper-branch coefficient is $10.1268$, again matching the leading constant as well as the exponent.

\begin{figure}[ht]
\centering
\begin{tikzpicture}[scale=2.0,>=Stealth]
\pgfmathdeclarefunction{dneutral}{1}{\pgfmathparse{0.8 + 0.14476636480382166*cos(#1) + 0.12866554*sin(#1) - 0.16790919*cos(2*#1) + 0.08267596*sin(2*#1) - 0.11676398*cos(3*#1)}}
\draw[thick] (0,0) circle (1);
\draw[thick] plot[smooth,domain=0:360,samples=240] ({(1+dneutral(\x))*cos(\x)},{(1+dneutral(\x))*sin(\x)});
\draw[-{Stealth[length=2mm]},densely dashed] (5:1.08) arc[start angle=5,end angle=345,radius=1.08];
\node at (-0.15,0.18) {$C$};
\node at (1.25,1.55) {$\partial\Omega$};
\node[fill=white,inner sep=1pt] at (-0.78,-0.73) {$\gamma$};
\draw[-{Stealth[length=1.6mm]}] (1.660093,0.000000) -- (0.992996,0.118149);
\draw[-{Stealth[length=1.6mm]}] (1.128064,1.953865) -- (0.427187,0.904163);
\draw[-{Stealth[length=1.6mm]}] (-1.371271,0.791704) -- (-0.858245,0.513240);
\draw[-{Stealth[length=1.6mm]}] (-0.827490,-1.433254) -- (-0.404025,-0.914748);
\draw[-{Stealth[length=1.6mm]}] (0.945037,-1.636853) -- (0.422830,-0.906209);
\fill (1,0) circle (.025) node[below left,xshift=-1pt] {$c_0$};
\fill (1.660093,0) circle (.025) node[right] {$\Phi(c_0)$};
\fill (0.992996,0.118149) circle (.022) node[above left,xshift=-2pt,yshift=1pt] {$F(c_0)$};
\draw[<->] (1.02,-0.10) -- (1.64,-0.10) node[midway,below] {$a_*$};
\end{tikzpicture}
\caption{TikZ rendering of the numerically tuned neutral circular-core geometry \eqref{eq:d-neutral}. The arrows are actual inward-normal return segments for five base points.}
\label{fig:neutral-geometry}
\end{figure}
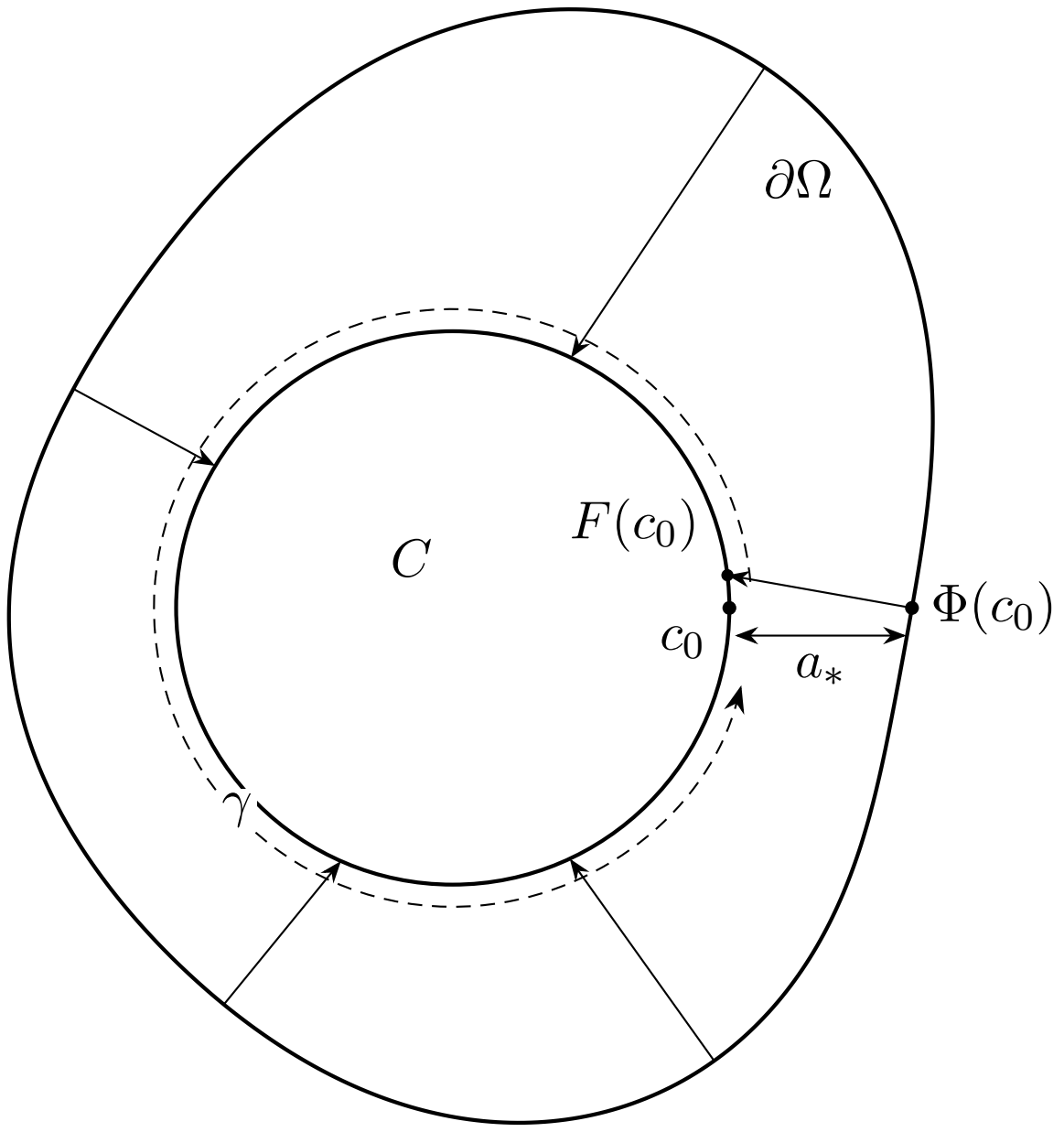

Figure~\ref{fig:neutral-geometry} displays the core, the Fourier outer boundary and several actual inward-normal return segments. The segment from $c_0$ to $\Phi(c_0)$ has length $a_*$ and the arrow on the core represents the closed loop used in the scalar transport.

Under the same perturbation \eqref{eq:q-noise}, the even-order neutral root splits into two branches. Table~\ref{tab:neutral-num} reports both anchor displacements and the global thickness error for the upper branch. The fitted slopes over the complete range $10^{-8}\le\rho\le3\times10^{-5}$ are
\[
0.49597\quad\text{for }|\hat a-a_*|,
\qquad
0.49528\quad\text{for }\|\hat d-d\|_\infty,
\]
in agreement with the theoretical exponent $1/2$.

\begin{table}[ht]
\centering
\caption{Tuned neutral geometric experiment.}
\label{tab:neutral-num}
\begin{tabular}{@{}cccc@{}}
\toprule
$\rho=\varepsilon_F$ & lower root error & upper root error & $\|\hat d-d\|_\infty$\\
\midrule
$10^{-8}$ & $1.01404\times10^{-3}$ & $1.01268\times10^{-3}$ & $1.09443\times10^{-3}$\\
$10^{-7}$ & $3.21135\times10^{-3}$ & $3.19773\times10^{-3}$ & $3.45493\times10^{-3}$\\
$10^{-6}$ & $1.02017\times10^{-2}$ & $1.00656\times10^{-2}$ & $1.08659\times10^{-2}$\\
$10^{-5}$ & $3.27262\times10^{-2}$ & $3.13661\times10^{-2}$ & $3.37753\times10^{-2}$\\
$3\times10^{-5}$ & $5.75431\times10^{-2}$ & $5.34719\times10^{-2}$ & $5.74398\times10^{-2}$\\
\bottomrule
\end{tabular}
\end{table}

\begin{figure}[ht]
\centering
\begin{tikzpicture}
\begin{loglogaxis}[
width=.72\textwidth,height=.43\textwidth,
xlabel={$\varepsilon_F$},ylabel={error},
legend pos=north west,grid=major]
\addplot[mark=*] coordinates {
(1e-8,1.01268265e-3) (3e-8,1.75315556e-3) (1e-7,3.19773182e-3)
(3e-7,5.53000539e-3) (1e-6,1.00655951e-2) (3e-6,1.73479383e-2)
(1e-5,3.13660733e-2) (3e-5,5.34718513e-2)};
\addlegendentry{$|\hat a-a_*|$}
\addplot[mark=square*] coordinates {
(1e-8,1.09443373e-3) (3e-8,1.89450629e-3) (1e-7,3.45492744e-3)
(3e-7,5.97305097e-3) (1e-6,1.08659473e-2) (3e-6,1.87108941e-2)
(1e-5,3.37752723e-2) (3e-5,5.74398145e-2)};
\addlegendentry{$\|\hat d-d\|_\infty$}
\addplot[dashed,domain=1e-8:3e-5,samples=2] {10*sqrt(x)};
\addlegendentry{reference slope $1/2$}
\end{loglogaxis}
\end{tikzpicture}
\caption{Tuned neutral return geometry: the anchor and globally reconstructed thickness both exhibit the square-root stability law.}
\end{figure}
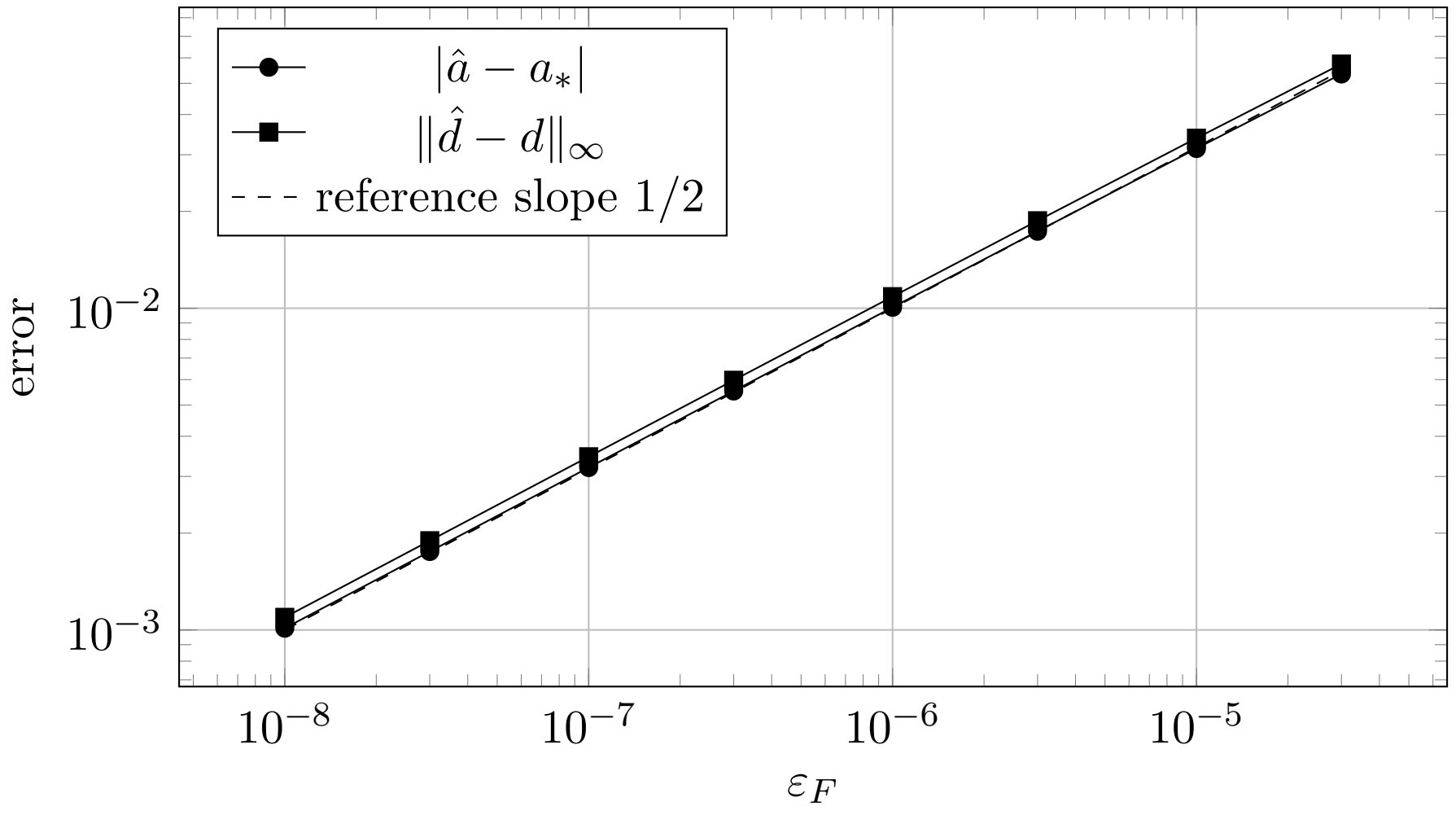

The observed separation $\min_\theta(\hat d(\theta)-\hat q_\rho(\theta))$ remains above $0.58$ at the smallest noise levels and above $0.63$ at $\rho=3\times10^{-5}$, so the experiment stays well inside the uniformly separated regime of the transport estimates. The numerics therefore validate not only the scalar root exponent but the full chain from perturbed return components to anchor reconstruction and global thickness propagation.

\section{Conclusion and open problems}
\label{sec:conclusion}

The paper develops a finite-order singular stability theory for the anchor parameter in the nonlinear holonomy formulation of the return-map inverse problem. The proved chain is
\[
\begin{aligned}
\text{return-map perturbation}
&\Longrightarrow \text{finite holonomy-jet perturbation},\\
&\Longrightarrow \text{finite isolation order}
\Longrightarrow \text{anchor stability},\\
&\Longrightarrow \text{multi-loop aggregation}
\Longrightarrow \text{global thickness stability}.
\end{aligned}
\]
At first order, the non-neutral regime is locally Lipschitz stable. At finite neutral order $k\ge2$, the anchor obeys a $1/k$-H\"older location estimate with parity-sensitive existence. Exact uniform separation is stable under sufficiently small return-data perturbations, so the singular fibre set enters the stability theory through an explicit quantitative conditioning margin. The exact second and third holonomy variations and the Bell-polynomial hierarchy identify the relevant coefficients as vertical derivatives of the geometric transport. For smooth finite-dimensional families, Malgrange preparation further shows that an order-$k$ holonomy singularity is locally equivalent, up to multiplication by a nonvanishing smooth factor, to a monic degree-$k$ polynomial in the anchor displacement; in the quadratic case Proposition~\ref{prop:malgrange-k2-coeff} identifies the first-order prepared coefficients and discriminant correctly, including the contribution of the cubic anchor jet to $c_1$. The finite-jet perturbation theorem makes the holonomy coefficients quantitatively observable from perturbed return data, while the multi-loop theorem combines several closed-loop signals without a statistical-independence assumption. Finally, the analytic-dependence theorem shows that infinite-order flatness for a fixed uniformly separated loop means local identity holonomy, rather than an arbitrary $C^\infty$ flat germ.

The numerical experiments support the theorem chain on genuine return geometries. A non-neutral circular example gives first-order anchor and global reconstruction rates, while a numerically tuned neutral Fourier geometry has $H'(a_*)$ at machine zero, $H''(a_*)\ne0$ and yields observed rates close to $1/2$ for both the anchor and the full reconstructed thickness. In both regimes the observed leading anchor constant agrees with the corresponding local normal-form prediction.

The principal remaining questions are geometric rather than scalar:
\begin{enumerate}[label=(\roman*)]
\item derive invariant formulas relating small-loop holonomy jets to the curvature and vertical curvature jet of the nonlinear inverse connection;
\item determine generic non-vanishing of those coefficients on genuine non-umbilic compact cores, including ellipsoids, without suppressing Levi-Civita connection terms;
\item characterize configurations for which all relevant loop holonomies are locally the identity and decide whether global admissibility or topology still isolates the anchor;
\item optimize loop families by combining isolation order, collective jet amplitude, loop length and transport conditioning;
\item extend the reconstruction and stability theory to lower regularity, for example $d\in C^{1,1}(M)$;
\item develop a range characterization for noisy candidate return maps that couples the transport residuals with geometric admissibility of the reconstructed normal graph.
\end{enumerate}

\end{document}